\documentclass{amsart}
\usepackage{amsfonts}
\begin{document}
\newtheorem{Theorem}{Theorem}[section]
\newtheorem{Remark}{Remark}[Theorem]
\newtheorem{Lemma}{Lemma}[Theorem]
\def\qedbox{\hbox{$\rlap{$\sqcap$}\sqcup$}}
\def\qed{\nobreak\hfill\penalty250 \hbox{}\nobreak\hfill\qedbox}
\def\e{\epsilon}
\def\diag{\mbox{diag}}
\def\dim{\mbox{dim}}
\def\arg{\mbox{arg}}
\def\tr{\mbox{tr}}
\def\e{\epsilon}
\def\ve{\varepsilon}
\def\RR{\mathbb{R}}
\def\CC{\mathbb{C}}
\def\HH{\mathbb{H}}
\def\FF{\mathbb{F}}
\def\NN{\mathbb{N}}
\def\II{\mathbb{1}}
\def\enabla{\stackrel{\varepsilon}{\nabla}}
\def\bnabla{\overline{\nabla}}
\def\cnabla{\stackrel{\circ}{\nabla}}
\def\lh{{\mathcal K}}
\def\lho{{\mathcal K}_0}
\def\L#1{\Lambda #1}
\def\lll{{\mathcal L}}
\def\O{{\mathcal O}}
\def\llh{{\mathcal L}_H/S^1}
\def\oo{{\mathcal O}_{p,q}(1)}
\def\ll{{\Lambda}}
\title{Examples of Self-dual, Einstein metrics  of $(2,2)$-signature}
\author{Novica Bla\v zi\' c and  Srdjan Vukmirovi\' c}
\address{Faculty of mathematics, University of Belgrade, Studenski trg 16, p.p.
550,
11 000 Belgrade, Yugoslavia}
\email{blazicn@matf.bg.ac.yu,  vsrdjan@matf.bg.ac.yu}
\subjclass{53C07, 53G10, 53C80, 53C15}
\keywords{Osserman manifold, indefinite metric}
\date{}
\begin{abstract}
  In this paper we construct a family of examples of self-dual Einstain metrics of neutral signature, which are  not Ricci flat, nor locally homogenous.
Curvature of these manifolds is studied in details.
These are  obtained by the   para-quaternionic reduction.
We compare our examples with the  orbifolds $\oo$ given by Galicki and Lawson,
for which some new  properties are also established.
Particularly, the  sign and  the pinching of their   sectional curvatures
are studied.
\end{abstract}
\maketitle
\section{Introduction}
A self-dual manifolds are important in the Riemannian geometry
 (of the positive definite  $(++++)$ signature), as well in the Kleinian geometry
(of the  neutral $(--++)$ signature).
Some examples of Einstein, self-dual metrics of neutral
signature are studied in the recent times~\cite{KM}, but non-trivial examples of
non-Ricci flat, Einstein self-dual manifolds were not known.
In the  positive definite case  a family of such examples
is constructed  in \cite{G} and \cite{GL}.
Motivated by these examples, here we construct
a family  $M_{p,q}$, $p,q\in \NN$, $(p,q)=1$, $p \neq q$  of Einstein, self-dual
 manifolds  of neutral signature, which are  even  not
locally homogenous.
As remarked in \cite{S}, the quaternionic reduction can be generalized to the pseudo-Riemannian
cathegory. Particularly, we apply the reduction procedure in the setting of
$C(1,1)$ Clifford algebra
(para-quaternionic numbers) instead of $C(2)$ Clifford algebra (quaternions).
This approach is appropriate  because of the  following reasons.
Curvature of the self-dual manifolds share some properties with
general quaternionic and para-quaternionic manifolds (see \cite{A, G, GL}).
The other reason is that the Clifford algebra $C(1, 1)$ provides
  a natural language for describing metrics of the neutral signatures.
It is interesting to compare this construction with example given in \cite{GL}
obtained by quaternionic reduction.
In the  both cases an algebraic submanifold $\lho$ of projective plane ($\HH P^2$ or $\tilde \HH P^2$) is constructed. In the quaternionic  (positive definite) case  $\lh = \lho$ is a differentiable, complete, submanifold
 but in the other case (para-quaternionic, indefinite)  $\lho$ is not a differentiable submanifold. Its singular set ${\mathcal S}$ lies on a real hyperquadric in $\tilde \HH P^2.$ Let $\lh = \lho \setminus{\mathcal S}.$
In  the both cases    the family of isometric actions $\phi _{p,q}, \enskip p, q \in \NN, \enskip (p,q) = 1$ of a group $G$
 ($S^1$ or $\RR$)  on five-dimensional submanifold $\lh$ is defined.
In the quaternionic case this is  the locally free action   and in
the para-quaternionic  case the actions $\phi _{p,q}$ are always
free. In the quaternionic case $\lh /G$ is a Riemannian  orbifold
and in para-quaternionic  case $\lh /G$ is a manifold with
uncomplete metric of  neutral signature. Their "singularity sets"
are  of different character.
Here is the brief outline to this paper.
In \S2, we shall discuss  the elementary  properties of the Clifford algebra
$\tilde \HH = C(1,1) = C(0,2)$ (para-quaternions) and the  projective plane which
it defines.
In \S3 we shall define the action of the group $G=\RR$  on the para-quaternionic projective plane
$\tilde{\HH}P^2$ and the algebraic submanifold  $\lh_0$
in that plane. Moreover, we shall study  the basic geometry of the orbit space $\lh_0/G$, equiped with the metric $g_{p,q}$ induced by the submersion  $\lh_0\to \lh_0/G$.
In \S4 we shall compute the sectional curvature and
 the Jacoby curvature operator of the metric $g_{p,q}$ in terms of
the para-quaternions. This shall give us  possibility to
 prove in \S5 that  $g_{p,q}$
is  the self-dual, Einstein metric (pointwise Osserman) which is not
even locally homogeneous.
We conclude  in \S6 by applying this ideas to study   in more details the  curvature
of the orbifolds $\oo$, defined by Galicki and Lawson (\cite{G,GL}). For some known results new proofs and generalizations are given.
\section{Preliminaries}
Denote by $ \tilde \HH = C(1,1) = C(0,2)$ the real Clifford algebra
 with unity $1$ of para-quaternionic numbers generated by elements
$i,$ $j,$ $k$ satisfying
$$
i^2 = -\e _1 = -1, \enskip j^2 = -\e _2 = 1, \enskip
k^2 = -\e _3 = 1, \enskip       ij=k=-ji.
$$
The conjugate $\bar q$ of para-quaternionic number $q := x + yi + zj + wk$
is defined
by
$$
\bar q :=  x - yi - zj - wk
$$
and real and imaginary part, respectively, by
$$
\Re q :=x, \enskip \Im q = yi + z j + w k.
$$
The square  norm of the  para-quaternion $q$
$$
|q|^2 := q \bar q= x^2 + y^2 -z^2 - w^2,
$$
is multiplicative,  i.e. satisfies
$$
|q_1 q_2|^2 = |q_1|^2 |q_2|^2 .
$$
Notice that  the pseudosphere  of signature $(2,1)$ of all unit para-quaternions
$$
S^{2,1} = \tilde \HH _1 := \{ q \in \HH \mid |q|^2 = 1 \}
$$
is a Lie group $\tilde \HH _1 = SU(1,1).$
On the   vector space $\tilde \HH ^3$ we have the pseudo-Riemannian metric of
signature $(6,6)$ given by
$$
g( u, v)  := \Re (\bar u_1v_1 + \bar u_2v_2 + \bar u_3v_3)
= \Re (\bar u\cdot v)
$$
for all $ u = (u_1, u_2, u_3), \enskip v = (v_1, v_2, v_3)$ in $\tilde  \HH ^3.$
Denote by $S^{6,5}$ pseudosphere of signature $(6,5)$  in $\tilde \HH ^3.$
The para-quaternionic projective plane $\tilde \HH P^2$
is defined as a set of equivalence classes
$$
\tilde \HH P^2 := \{ [u] \mid u \in S^{6,5} \}
$$
where the  equivalence $\sim$ is defined  by
$$
u = (u_0, u_1, u_2 ) \sim (u_0h, u_1h, u_2h) = uh, \enskip h
\in \tilde \HH _1.
$$
Denote by $\pi : S^{6,5} \to \tilde \HH P^2$ the natural projection $\pi (u):=[u].$
The vertical subspace $ {\mathcal V}_u\subset T_uS^{6,5} \subset T_u\tilde \HH ^3 \cong \tilde \HH ^3$ of the  submerison $\pi$
in the  point $u \in S^{6,5}$ is generated by the vectors $ui, uj, uk$ i.e.
$$
{\mathcal V}_u = \RR \langle ui, uj, uk \rangle .
$$
Moreover, the orthogonal decomposition with respect to the  metric $g$
$$
T_uS^{6,5} = {\mathcal V}_u\oplus {\mathcal H}_u
$$
onto vertical and horizontal subspace is $\tilde \HH _1$  invariant.
We usually fix a horizontal local  section $u \subset S^{6,5}$ over $\tilde \HH P^2 .$ This allows us to indentify
$$
T_{[u]}\tilde \HH P^2 \cong {\mathcal H}_u \enskip \mbox{ by } \enskip T_{[u]}\tilde \HH P^2
\ni X \leftrightarrow \tilde X \in {\mathcal H}_u \subset \tilde \HH ^3,
$$
where $\tilde X$ is the unique horizontal lift of the vector $X.$
Right multiplication of vectors form $\tilde \HH ^3$ by $i, j, k$ preserves  horizontal subspace ${\mathcal H}_u$ so we can locally define the endomorphisms
$$
J_1(X) := \pi _*(\tilde Xi), \enskip J_2(X) := \pi _*(\tilde Xj), \enskip J_3(X) := \pi _*(\tilde Xk)
$$
where the  vectors $X$ and $\tilde X$ are related as above. They satisfy the same relations as $i,j$ and $k$, i.e.
\begin{equation}
J_\alpha ^2 = -\e_\alpha , \enskip J_\alpha J_\beta  = -\e_\gamma J_\gamma ,
\label{com}
\end{equation}
where $(\alpha, \beta, \gamma)$ is any cyclic permutation of $(1, 2, 3).$
Although the  definition of $J_1, J_2$ and $J_3$ depends on the section $u$ over $\tilde \HH P^2$ we have globally defined  subbundle
$\RR \langle J_1, J_2, J_3 \rangle $ of $End(T\tilde \HH P^2)$ which is a pseudo-Riemannian analogue of quaternionic K\" ahler structure on quaternionic projective plane $\HH P^2.$
The pseudo-Riemannian metric $g$
on $S^{6,5}$ induces  the pseudo-Riemannian metric
$$
\langle X, Y \rangle = g(\tilde X, \tilde Y)
$$
of signature $(4,4)$ on
$\tilde \HH P^2$ by the submersion.
One can easily check that the endomorphism $J_1$ is an isometry and $J_2, J_3$ are anti-isometries.
We denote the Levi-Civita connection on $S^{6,5}$  by $\enabla$
and by $\bnabla$ the induced connection  on  $\tilde \HH P^2$ with a
curvature tensor $\bar R$ and of the  constant para-quaternionic sectional
curvature $\bar K$ (see~\cite{B}).
\section{The orbit space $\lh /G$ and its metric}
Let $p, q \in \NN, (p, q)=1, p \neq q$ be relatively prime natural numbers.
We define an isometric action $\phi$ of a group
$G := \{ e^{jt} \mid t \in \RR \} \cong (\RR, +)$
on $\tilde \HH P^2$ by
$$
\phi _t (u_0, u_1, u_2) := (e^{jqt}u_0, e^{jpt}u_1, e^{jpt}u_2),
$$
where $ e^{jt} := \cosh t + j\sinh t,\enskip t\in \RR .$
The action preserves para-quaternionic structure on $\tilde \HH P^2.$
The induced Killing vector field at the point $[u]\in \tilde \HH P^2$ is
$$
 V_u := \pi_*(j(qu_0, pu_1, pu_2)).
$$
Consider an algebraic submanifold $\lho$ of $\tilde \HH P^2$ defined
 by equation
$$
\lho := \{ [u_0, u_1, u_2] \in \tilde \HH P^2 \mid
q\bar u_0ju_0 + p
\bar u_1ju_1 + p\bar u_2ju_2 = 0\}.
$$
Notice that its real  codimension in $\tilde \HH P^2$ is  $3$ and  $dim_\RR \lho = 5.$
\begin{Lemma}
\label{free}
The action $\phi _t$ of the  group $G$ preserves the submanifold $\lho$.
It acts freely on $\lho.$
\end{Lemma}
Proof: The manifold $\lho$  is invariant under the  the action of $G$ because of the definition of the definition of $\phi_{p,q}$.
Suppose that the   action of $G$ is not free, i.e. that for some $h \in \HH _1,$ $t \in \RR$, $t\neq 0$ and some $u \in \lho \subset S^{6,5}$ we have
$$
(e^{jqt}u_0, e^{jpt}u_1, e^{jpt}u_2) = \phi _t (u) = uh = (u_0h, u_1h, u_2h).
$$
Using the equation of the  submanifold $\lho$ we obtain
\begin{eqnarray*}
0 & = & q\bar u_0ju_0h + p\bar u_1ju_1h + p\bar u_2ju_2h = \\
& = &
q \bar u_0je^{jqt}u_0 + p\bar u_1je^{jpt}u_1+ p\bar u_1je^{jpt}u_2 =\\
& = &
 |u_0|^2 (q\sinh qt-p\sinh pt) + p\sinh pt
+ q(\cosh qt - \cosh pt)\bar u_0 ju_0 .
\end{eqnarray*}
Hence, we have $q(\cosh qt - \cosh pt)\bar u_0 ju_0  \in \RR$.
If $\cosh qt = \cosh pt$ then, since $t \neq 0,$ we
have $p = \pm q,$ what is a contradiction.
Otherwise, we have  $\bar u_0 ju_0 = a \in \RR.$
If $|u_0|^2 \neq 0$ then $j = a|u_0|^{-2}\in \RR.$
It is a contradiction again.
If $|u_0|^2 = 0$ then $a \neq 0,$ since $ p\sinh pt \neq 0,$ so we have
$$
0 = |u_0|^2|j|^2|u_0|^2 = |a|^2 \neq 0
$$
and the action of $G$ is free.
\qed
Since the action $\phi _t$ of the group $G$ preserves
the algebraic submanifold $\lho$,
the vector $V_u$ is tangent to $\lho$ at the point $[u] \in \lho.$
One can check that the vectors
$$
J_1V_u, \enskip   J_2V_u, \enskip   J_3V_u \in T_{[u]}\tilde \HH P^2
$$
are normal to  the submanifold $\lho$ at the  point $[u].$  In
points $[u] \in \lho$ such that  $|V_u|^2 \neq 0$
 vectors $V_u, J_1V_u, J_2V_u, J_3V_u$  are linearly independent.
Hence, in order to proceed with calculations we restrict to the
subset  $\lh$ of $\lho$ consisting of such points, i.e.
$$
\lh := \{ [u]=[u_0, u_1, u_2] \in \lho \mid n^2 :=|V_u|^2 = q^2 |u_0|^2 + p^2 |u_1|^2
+ p^2 |u_2|^2 \neq 0\}.
$$
On $\lh$ we have the  metric $\langle  \cdot, \cdot \rangle$ induced from $\tilde \HH P^2$ and the
connection, curvature and sectional curvature  which  we denote by
$\cnabla, \stackrel{\circ}{R},\stackrel{\circ}{K}, $ respectively.
The following lemma is immediate consequence of Lemma \ref{free}.
\begin{Lemma}
Set $\lh $ is a differentiable manifold of real dimension five.
The isometric action $\phi _t$ of the group $G$ preserves $\lh$ and it is free on $\lh.$
Then $\lh /G$ is a pseudo-Riemannian manifold of real signature $(2,2).$
\qed
\end{Lemma}
Denote by $g_{p,q}$ the metric induced on $\lh /G$ by Riemannian submersion
$\xi : \lh \to \lh /G.$ The construction  is
represented by the following diagram.

\vspace{2em} \noindent {\small
\begin{tabular}{lllllllllll}
  &   &\hfil$\tilde \HH_1=S^{2,1}$\hfil & $\to$ & $G=\RR$ & $\to$ & $G=\RR$ &
  $\to$  & $G=\RR$         &     &    \\
       &     &   \hfil$\downarrow$\hfill &     &\hfil$\downarrow\phi_{p,q}$\hfil &
&\hfil$\downarrow\phi_{p,q}$\hfil &
&\hfil$\downarrow\phi_{p,q}$\hfil &  &      \\
$\tilde\HH^3$ & $\supset$ &\hfil $S^{6,5}$ \hfil
& $\stackrel{\pi}{\rightarrow}$ &
$\tilde \HH P^2$ &  $\supset $ & $\lho$   &   $\supset $
& $\lh=\lho\setminus {\mathcal S}$
& $\stackrel{\xi}{\rightarrow}$ &
$M_{p,q}=\lh/G$      \\
\end{tabular}
} 
\vspace{2em}

\noindent It is not difficult to check that it can be expressed in
terms of global coordinates $u \in S^{6,5}$ in the following way
$$
(\xi \circ \pi)^*(g_{p,q}) = d\bar{u}\otimes_\RR du
+ (d\bar{u}\cdot u )\otimes_\RR  (\bar u \cdot du)
 +
\frac{(d\bar{u}\cdot V_u) \otimes_\RR
( \overline{V_u}\cdot du)}{|V_u|^2 } .
$$
Notice, that the first two terms determine the metric of the
para-quaternionic projective plane.
The endomorphisms $J_i, \enskip i=1,2,3$ of $T_{[u]}\tilde \HH P^n$ induce  the endomorphisms $\tilde J_i, \enskip i=1,2,3$ of $T_{\xi ([u])}(\lh /G)$, satisfying the relations (\ref{com}),
 by
$$
\tilde J_i(\tilde X, \tilde Y ) = J_i(X,Y), \enskip \tilde X, \tilde Y \in T_{\xi ([u])}(\lh /G),
$$
where $X, Y$ are horizontal lifts of vectors $\tilde X, \tilde Y$.
Moreover, $\tilde J_1$ is the isometry and $\tilde J_2, \tilde J_3$ are the anti-isometries  with respect to the metric $g_{p,q}.$
Denote the induced connection, curvature and the
 sectional curvature on $\lh /G$ by $ \nabla,$ $R$ and $K,$ respectively.
\section{The curvature of $\lh /G$}
In order to study the local geometry of $\lh /G$ we  will  compute
 its curvature tensor in this section.
The sectional curvature of the quaternionic and  the para-quaternionic
projective plane can be expressed as
\begin{equation}
\bar K(X,Y)=1+3|Pr_{\langle Y,J_1Y,J_2Y,J_3Y\rangle}X|^2/\epsilon(X),
\label{sc}
\end{equation}
for arbitrary orthonormal tangent vectors $X$ and $Y,$ where $\epsilon(X)$ denotes the sign of square  length of the  vector $X$ (see~\cite{B}).
Now, we will compute the connection $\cnabla$ and
the curvature $\stackrel{\circ}{R}$ of the  submanifold  $\lho $ of $\tilde  \HH P^2.$
Let $i : \lh \to  \tilde \HH P^2$ denotes the  inclusion.
For  $[u] \in \lh$ we will denote the vectors $X \in T_{[u]}\lh$ and $i_*X \in T_{[u]}\tilde \HH P^2$ by the same letter $X.$
The following relation holds
\begin{equation}
\label{be}
\bnabla_XY = \cnabla_XY + B(X, Y), \enskip X, Y \in T_{[u]}\lh ,
\end{equation}
where $B(X, Y)$ is the  second fundamental form at  a point $[u] \in \lh .$
Since $B(X, Y) \in \left( T_{[u]}\lh \right) ^\perp
= \RR \langle J_1V_u, J_2V_u, J_3V_u \rangle $
we have
$$
\bnabla_XY - \cnabla_XY =  B(X, Y)
= \sum_{i=1}^3\alpha_i J_iV_u
$$
where
\begin{equation}
\label{abc}
\alpha_i = \frac{\e _i}{n^2} \langle \bnabla_XY, J_iV_u \rangle ,\quad
i=1,2,3.
\end{equation}
For  $i=1,2,3$ we have
\begin{equation}
\label{nabla}
 \langle \bnabla_XY, J_iV_u \rangle  = - \langle Y, \bnabla_XJ_iV_u \rangle  =
- \langle Y, J_i\L{X} \rangle ,
\end{equation}
where  the skew-symmetric operator $\Lambda  : T_{[u]}\tilde \HH P^2 \to T_{[u]}\tilde \HH P^2$ is defined by
$$
\Lambda X :=\pi _*(j(qx_0, px_1, px_2))
$$
and $(x_0, x_1, x_2)$ is horizontal lift of $X \in T_{[u]}\tilde \HH P^2$ to $T_uS^{6,5}.$
Using equations (\ref{be}), (\ref{abc}) and (\ref{nabla})
we finally obtain
$$
B(X, Y) = -\frac{1}{n^2} \sum_{i=1}^{3} \e_i  \langle Y, J_i\L{X} \rangle  J_iV_u.
$$
Using the relation
\begin{eqnarray*}
\stackrel{\circ}{K}(X, Y)  & = & \bar K(X, Y) -\frac{1}{Q}(|B(X, Y)|^2 -  \langle B(X, X), B(Y, Y) \rangle ), \\
 Q  &:=&  Q(X, Y) = |X|^2|Y|^2 -  \langle X, Y \rangle ^2,
\end{eqnarray*}
between second fundamental form and sectional curvature we obtain
$$
\stackrel{\circ}{K}(X, Y)  = \bar K(X, Y) +
\frac{1}{n^2Q} \sum _{i=1}^3 \e_i \left( - \langle Y, J_i\L{X} \rangle ^2  +  \langle X, J_i\L{X} \rangle  \langle Y, J_i\L{Y} \rangle  \right),
$$
where the sectional curvature $\bar K$ is constant.
As the  second step, we are going to compute the
 sectional curvature $K$ on $\lh  /G$ using O'Neill's  formula for submersion
\begin{equation}
\label{k(xy)}
K(\tilde X, \tilde Y) = \stackrel{\circ}{K}(X, Y)  + \frac{3}{4}\frac{ \langle v[X, Y], v[X, Y] \rangle }{Q(X, Y)},\enskip \tilde X, \tilde Y \in T_{\xi ([u])}(\lh /G)
\end{equation}
where $v[X, Y]$ denotes the vertical component of the  commutator $[X, Y]$ and $X, Y$ are the unique horizontal lifts   of  the vectors $\tilde X, \tilde Y$ to $T_{[u]}\lh .$
Since the  vertical space of the submersion is generated by the vector $V_u$
we have
\begin{equation*}
\begin{aligned}
v[X, Y] & =  \frac{1}{n^2} \langle [X, Y], V_u \rangle  V_u
= \frac{1}{n^2} ( \langle \cnabla_XY, V_u \rangle
-  \langle  \cnabla_YX, V_u \rangle )V_u = \\
     & =  \frac{1}{n^2}(- \langle Y, \L{X} \rangle  +  \langle X, \L{Y} \rangle)V_u
=\frac{2}{n^2}  \langle X, \L{Y} \rangle V_u .
\end{aligned}
\end{equation*}
Using (\ref{k(xy)}) the  sectional curvature of $\lh /G$ is
\begin{equation*}
\begin{aligned}
K(\tilde X, \tilde Y) & =  \bar K(X, Y) +\\
& +   \frac{1}{n^2 Q}\left[ \sum_{i=1}^3 \e _i
\left(   \langle X, J_i\L{X} \rangle  \langle Y, J_i\L{Y} \rangle
-  \langle Y, J_i\L{X} \rangle ^2\right)  + 3 \langle X,\L{Y} \rangle ^2 \right]
\end{aligned}
\end{equation*}
and $X, Y$ are related to $\tilde X, \tilde Y$ as above.
Then, the Jacoby operator is
\begin{equation*}
\begin{aligned}
K_{\tilde X}(\tilde Y) &
:=  R(\tilde X, \tilde Y)\tilde X   = h(\bar R(X, Y)X)  +
 \frac{1}{n}
\left[
\sum _{i=1}^3\e _i \left( - \langle X, J_i\L{X} \rangle
h(J_i\L{Y}) +\right. \right.  \\
&\left. \left.   +  \langle X, J_i\L{Y} \rangle h(J_i\L{X}) \right) - 3 \langle Y, \L{X}
\rangle h(\L{X}) \right],
\end{aligned}
\end{equation*}
where $h(X)$ denotes the  horizontal part of the  vector $X \in T_{[u]}\tilde \HH P^2$ , i.e. the orthogonal  projection of the  vector $X$ onto the  space $\RR \langle V_u, J_1V_u,
J_2V_u, J_3V_u \rangle ^\perp .$
\section{The local geometry of $\lh /G$}
In this section we shall use the curvature of $\lh /G$
to describe  its    interesting  geometric properties.
First, we recall the definition of the  Osserman manifold $(M, g).$
Let $m \in M$ and $X \in T_mM$ be a vector such that $|X|^2 = 1.$ Jacoby operator $K_X$ at the   point $m \in M$ and in the direction $X$  is defined by
$$
K_X(Y) := R(X,Y)X,
$$
where $R$ is the curvature tensor of $(M, g).$
Jacoby operator is a self-adjoint operator on $T_mM.$
Manifold $M$ is called  pointwise Osserman if the  Jordan form of the Jacoby operator $K_X$ does not depend on the unit direction $X$
and (globally) Osserman if the  Jordan form  of $K_X$ does not depend both on the  point $m \in M$ and the  unit direction $X \in T_mM$.
\begin{Lemma}
\label{l51}
The manifold $\lh /G$ is pointwise Osserman of the  neutral signature.
\end{Lemma}
Proof: Let $X \in T_{[u]}\lh$ denote the horizontal lift of a vector $\tilde X \in T_{\xi ([u])}(\lh /G).$
For any unit vector $\tilde X \in T_{\xi ([u])}(\lh /G)$ (i.e. with square
  norm equals $1$) the basis $(\tilde X, \tilde J_1\tilde X,\tilde J_2\tilde X, \tilde J_3\tilde X)$ is
 pseudo-orthonormal.
The vector $\tilde X$ is the eigenvector  corresponding to the  eigenvalue $0$ of the
self-adjoint operator  $K_{\tilde X}$. Thus, we are interested only in the  restriction
of the operator $K_{\tilde X}$ to the orthogonal complement
$\RR  \langle \tilde J_1\tilde X,\tilde J_2\tilde X, \tilde J_3\tilde X \rangle $ of $\tilde X.$
The restriction of $K_{\tilde X}$ in the  basis $(\tilde J_1\tilde X,\tilde J_2\tilde X, \tilde J_3\tilde X)$ is
represented by the following matrix
$$
K = \frac{2}{n^2}\left(
\begin{array}{ccc}
2a^2 + b^2 + c^2 -\frac{\bar{c}n^2}{2}   &    -3ab    &      -3ac    \\
3ab    &      -a^2 - 2b^2 + c^2 -\frac{\bar{c}n^2}{2}    & -3bc\\
3ac & -3bc & -a^2 + b^2 - 2c^2 -\frac{\bar{c}n^2}{2}
\end{array}
\right),
$$
where $a, b, c$ are the  coordinates of the horizontal
component $h(\L{X})$ of the vector $\L{X}$ in the  basis $(J_1X,J_2X, J_3X)$, i.e.
$$
a:=  \langle X, J_1\L{X} \rangle , \quad b:=  \langle X, J_2\L{X} \rangle ,
\quad c:=  \langle X, J_3\L{X} \rangle ,
$$
and $\bar{c}$  is the constant  para-quaternionic sectional curvature of
the corresponding projective plane.
By a  direct computation  we find that the  eigenvalues of the  matrix $K$ are
$$
\lambda _1 = -\frac{2}{n^2}(-a^2 + b^2 + c^2) +\bar{c}= \lambda _2,
\quad \lambda _3 = \frac{4}{n^2}(-a^2+b^2+c^2)+\bar{c} .
$$
To show that the  space $\lh /G$ is the  pointwise Osserman we have to show that the
eigenvalues $\lambda _i$ are independent of the chosen unit direction  $\tilde X,$ i.e. of its horizontal unit lift $ X \in T_{[u]}\lh$.
Notice that  $-a^2 + b^2 + c^2  = |h(\L{X})|^2$.
Let $e \in T_{[u]}\lh ,$ $|e|^2 = 1$ be a fixed unit
vector and $(e, J_1e, J_2e, J_3e)$
the corresponding pseudo-orthonormal basis.
Any horizontal unit vector $X$ can be written in the form
$$
X = X_0e + X_1J_1e +  X_2J_2e+  X_3J_3e, \enskip X_0^2+X_1^2-X_2^2-X_3^2=1.
$$
Now one can check directly that
$$
-a^2 + b^2 + c^2  = |h(\L{X})|^2 = |h(\L{e})|^2 = const
$$
and hence the eigenvalues of the Jacobi operator $K_{\tilde X}$ are constant
at a  given point.
Moreover, one can check that the Jacobi operator is  diagonalizable what completes the proof of the lemma.
\qed
\begin{Theorem}
The manifold $\lh /G$ is
  Einstein,  self-dual and not locally homogenous of the  signature $(--++)$.
\label{t51}
\end{Theorem}
Proof:
In \cite{ABBR},
the self-dual, Einstein manifold  of neutral signature
are characterized
as the pointwise-Osserman manifolds.
Then Lemma~\ref{l51} implies  that $\lh /G$ is Einstain and self-dual.
A pointwise-Osserman, locally homogenous manifold is
Osserman, i.e. the  Jordan normal  form of its Jacobi operator
is also independent of the point.
But for  our examples, we can  check that  the eigenvalues of the Jacoby operator, $\lambda_1$, $\lambda_2$ and $\lambda_3$ are
 are not constant on $\lh /G.$
After some computations, where  a few  unexpected cancelations
happened, we obtain
\begin{equation}
\lambda_1=\lambda_2=-\frac{2p^4q^2}{n^6}-4,\quad
\lambda_3=\frac{4p^4q^2}{n^6}-4. \label{51}
\end{equation}
Clearly, $n^2$ is not a  constant  along $\lh /G$ ($\lh /G$ is not
globally Osserman). Thus, the manifold $\lh /G$ is not locally
homogeneous, and hence
 not locally symmetric.
\qed

\noindent
{\bf Remark:} It is interesting that there exist a three-dimensional submanifold of $\lh /G$ along which the eigenvalues of the Jacoby operator are constant.

\noindent {\bf Remark: } In the case $p=1=q$ there are no singular
points, so  $\lh = \lho.$ Since the action of the group $G$ on
$\lh$ is free the resulting space $\lh /G$ is a complete manifold
of neutral signature. One can check that vector $\L{X},\enskip
X\in T_{[u]}\lh $ is horizontal and
 $\bar c = 4,$ so   $\lambda _1 = 2 = \lambda _2$, $\lambda _3 = 8$, i.e.  the eigenvalues of the  Jacoby operator are constant over $\lh /G.$
This implies that $M_{1,1}=\lh /G$ is Einstein, self-dual, globally Osserman neutral manifold. Moreover, it is  isometric to para-complex projective plane, where  para-complex stucture is induced by the  left multiplication by the generator $j$ (see~\cite{BBR}).
\section{Study of the  Galicki and Lawson's example}
We are coming back to the starting point
for our work, the compact Riemannian  orbifolds
$$
\O_{q,p}=\llh,
$$
constructed and studied  by Galicki and Lawson~\cite{GL}.
In this section we will discuss some additional  properties
of $\oo$ concerning to its local geometry and curvature.
Particularly, we will see that the  orbifolds $\oo$ are not locally homogeneous for any  $p,q$
and the orbifolds are of  positive sectional curvature  for some
 values of $p,q$. Also, an estimate of the pinching constant, ratio of the minimal and maximal
sectional curvature,  is given.
In \cite{GL} was shown that  $\O_{q,p}(n-1)$
is not locally symmetric for $0<q/p<1$.
Here, for $n=2$  we will give another   proof  of this
 result based  on the study of its
curvature  operators (as in Section \S5).
This approach provides an  oportunity  to see that $\oo$ is not even locally homogeneous.
\begin{Theorem}[\cite{GL}]
The orbifold $\O_{q,p}=\llh$
 is an Einstein, self-dual manifolds
which   is not  locally  homogeneous for any    $p,q$, $(p,q)\neq 1$.
\end{Theorem}
Proof:
The  proof is  based on the characterization  of
the Einstein, self-dual manifolds as manifolds with
constant   eigenvalues of the Jacobi operators
(pointwise Osserman manifolds), obtained by Vanhecke
and Sekigawa~\cite{SV}.
Repeating  the  long and nice
computations from the neutral signature case
(under some minor  modifications), we express the eigenvalues of the
Jacobi operator, $K_{\tilde X}:Y\to R(\tilde X,Y)\tilde X$,  in very simply way
\begin{equation}
\lambda_1=\lambda_2=\frac{2p^4q^2}{|V_u|^6}-4,\quad
\lambda_3=-\frac{4p^4q^2}{|V_u|^6}-4, \label{61}
\end{equation}
at  the orbit  determined by $[u]\in \lll_H$
 and the arbitrary unit tangent
vector $\tilde X$.
In this computations  operator
$\ll:T_{[u]}\tilde \HH P^2\to T_{[u]}\tilde \HH P^2$
plays an important role.
It is important to notice  that we express
the eigenvalues in the terms of the quaternionic notation.
By the argument as in the proof of Theorem~\ref{t51}
proof is completed.
\qed
Now, we will study the cases when the sectional
curvature of the
orbifold $\oo$ is   positive.
\begin{Lemma}
Let $\tilde X$ and $\tilde Y$  be the unit orthogonal tangent vectors.
For the sectional curvature $K(\tilde X,\tilde Y)$ of the
orbifold
$\oo$,  for $p\leq q$, holds
\begin{equation}
4-\frac{2q^2}{p^2}\leq K(\tilde X,\tilde Y)\leq 4+\frac{4q^2}{p^2},
\label{62}
\end{equation}
and for $p\geq q$
\begin{equation}
4-\frac{2p^4}{q^4}\leq K(\tilde X,\tilde Y)\leq 4+\frac{4p^4}{q^4}.
\label{63}
\end{equation}
\end{Lemma}
Proof:
For $p\leq q$ we have
$p^2\leq |\ll u|^2\leq q^2$
and
$$
-\lambda_1 \geq 4-\frac{2q^2}{p^2}, \quad
-\lambda_3\leq 4+\frac{4q^2}{p^2}.
$$
Since $-\lambda_1\leq K(x,y)\leq -\lambda_3$,
the relation (\ref{62}) holds.
Similarly,
for $p\geq q$ we have
$q^2\leq |\ll u|^2\leq p^2$
and
$$
-\lambda_1 \geq 4-\frac{2p^4}{q^4}, \quad
-\lambda_3\leq 4+\frac{4p^4}{q^4}.
$$
As in the previous case this  completes  the proof
the Lemma.
\qed

\noindent The direct consequences are the following results.
\begin{Theorem}
The orbifold $\oo$ is of the positive sectional curvature
for $ p^2 < \sqrt{2}q^2 < 2\sqrt{2} p^2$. \qed
\end{Theorem}
\begin{Lemma}
For $ p^2 < \sqrt{2}q^2 < 2\sqrt{2} p^2$, let the sectional
curvature $K(\tilde X, \tilde Y)$   of the orbifold $\oo$ is
$k$-pinched, $0<k<1$. Then, for $p\leq q$, in arbitrary point
$$
\frac{1}{4}-\frac{3}{4}\frac{q^2-p^2}{p^2+q^2}\leq  k \leq
\frac{1}{4}+\frac{3}{4}\frac{q^4-p^4}{p^4+q^4},
$$
and, for $p\geq q$
$$
\frac{1}{4}-\frac{3}{4}\frac{p^4-q^4}{p^4+q^4}\leq  k \leq
\frac{1}{4}+\frac{3}{4}\frac{p^2-q^2}{p^2+q^2}.
$$
For $p\rightarrow q$,  $k\rightarrow 1/4$.
Moreover,  $k=1/4$ on $\oo$ if and only if $p=q=1$.
\qed
\end{Lemma}
\noindent
{\bf Remark: }
Orbifold $\O_{1,1}(1)$ is globally Osserman and  of constant  holomorphic sectional
curvature with respect to  the complex structure
induced by the   left  multiplication by $j$.  (see \cite{C}).

\end{document}